\documentclass{amsart}
\usepackage{graphics,verbatim}	
\input epsf         

\newtheorem{proposition}{Proposition}[section]
\newtheorem{theorem}[proposition]{Theorem}

\newtheorem{lemma}[proposition]{Lemma}

\newtheorem{remark}[proposition]{Remark}
\newtheorem{example}[proposition]{Example}

\newtheorem{observation}[proposition]{Observation}

\newtheorem{prop}[proposition]{Proposition}
\newtheorem{cor}[proposition]{Corollary}

\numberwithin{equation}{section}

\newcommand{\Sh}{{\mathcal Sh}}
\newcommand{\rank}{\operatorname{rank}}

\newcommand{\Cg}{\mbox{{\rm Cg}}}
\newcommand{\Con}{\mbox{{\rm Con}}}
\newcommand{\Irr}{\mbox{{\rm Irr}}}
\newcommand{\cov}{\mathrm{cov}}
\newcommand{\covers}{{\,\,\,\cdot\!\!\!\! >\,\,}}
\newcommand{\covered}{{\,\,<\!\!\!\!\cdot\,\,\,}}
\newcommand{\set}[1]{{\left\lbrace #1 \right\rbrace}}
\newcommand{\pidown}{\pi_\downarrow}
\newcommand{\piup}{\pi^\uparrow}
\newcommand{\br}[1]{\langle #1 \rangle}
\newcommand{\A}{{\mathcal A}}
\newcommand{\join}{\vee}
\newcommand{\meet}{\wedge}
\newcommand{\leftq}[2]{\!\!\phantom{.}^{#1} {#2}}
\newcommand{\closeleftq}[2]{\!\!\phantom{.}^{#1}\! {#2}}

\begin{document}
\title{Sortable elements and Cambrian lattices}

\author{Nathan Reading}
\address{
Mathematics Department\\
    University of Michigan\\
    Ann Arbor, MI 48109-1043\\
USA}
\thanks{The author was partially supported by NSF grants DMS-0202430 and DMS-0502170.}
\email{nreading@umich.edu}
\urladdr{http://www.math.lsa.umich.edu/\textasciitilde nreading/}
\subjclass[2000]{20F55, 06B10}

\begin{abstract}
We show that the Coxeter-sortable elements in a finite Coxeter group $W$ are the minimal congruence-class representatives of a lattice congruence of the weak order on $W\!$.\@
We identify this congruence as the Cambrian congruence on $W,$ so that the Cambrian lattice is the weak order on Coxeter-sortable elements. 
These results exhibit $W$-Catalan combinatorics arising in the context of the lattice theory of the weak order on $W\!$.\@
\end{abstract} 

\maketitle

\setcounter{tocdepth}{1}
\tableofcontents

\section{Introduction}
\label{intro}

The weak order on a finite Coxeter group is a lattice~\cite{orderings} which encodes much of the combinatorics and geometry of the Coxeter group.
The weak order has been studied in the special case of the permutation lattice and in the broader generality of the poset of regions of a simplicial hyperplane arrangement.
With varying levels of generality, many lattice and order properties of this lattice have been determined.  
(See, for example, references in~\cite{bounded}, \cite{Du-Ch}, \cite{Markowsky}, \cite{hyperplane} and~\cite{hplanedim}.)  

This paper continues a program, begun in~\cite{congruence}, of applying lattice theory to gain new insights into the combinatorics and geometry of Coxeter groups.
Specifically, we solidify the connection, first explored in~\cite{cambrian}, between the lattice theory of the weak order and the combinatorics of the {\em \mbox{$W$-Catalan} numbers}.
These numbers count, among other things, the vertices of the (simple) generalized associahedron (a polytope which encodes the underlying structure of cluster algebras of finite type~\cite{ga,CA2}), the $W$-noncrossing partitions (which provide an approach~\cite{Bessis,BW} to the geometric group theory of the Artin group associated to~$W$) and the sortable elements~\cite{sortable} of~$W$ (which we discuss below).

In~\cite{cambrian}, the {\em Cambrian lattices} were defined as lattice quotients of the weak order on~$W$ modulo certain congruences, identified as the join (in the lattice of congruences of the weak order) of a small list of join-irreducible congruences.
Any lattice quotient of the weak order defines~\cite{con_app} a complete fan which coarsens the fan defined by the reflecting hyperplanes of $W,$ and in~\cite{cambrian} it was conjectured that the fan associated to a Cambrian lattice is combinatorially isomorphic to the normal fan of the corresponding generalized associahedron.
In particular, each Cambrian lattice was conjectured to have cardinality equal to the \mbox{$W$-Catalan} number.
These conjectures were proved for two infinite families of finite Coxeter groups ($A_n$ and~$B_n$).

The definition of {\em Coxeter-sortable elements} (or simply {\em sortable elements}) of~$W$ was inspired by the effort to better understand Cambrian lattices.
Sortable elements were introduced in~\cite{sortable} and used to give a bijective proof that $W$-noncrossing partitions are equinumerous with vertices of the generalized associahedron.
In this paper, we make explicit the essential connection between sortable elements and Cambrian lattices, proving in particular that the elements of the Cambrian lattice are counted by the \mbox{$W$-Catalan} number.
The conjecture from~\cite{cambrian} on the combinatorial isomorphism between Cambrian fans and cluster fans is proven in~\cite{camb_fan}.

The construction (see Section~\ref{sort sec}) of sortable elements involves the choice of a Coxeter element~$c$ of $W\!$.\@
For each~$c$, the corresponding sortable elements are called {\em \mbox{$c$-sortable}}.
The first main result of this paper is the following theorem.

\begin{theorem}
\label{cong}
Let~$c$ be a Coxeter element of a finite Coxeter group~$W\!$.\@
There exists a lattice congruence $\Theta_c$ of the weak order on~$W$ such that the bottom elements of the congruence classes of $\Theta_c$ are exactly the \mbox{$c$-sortable} elements.  In particular, the weak order on \mbox{$c$-sortable} elements is a lattice quotient of the weak order on all of~$W\!$.\,
\end{theorem}
The phrase ``bottom elements'' in Theorem~\ref{cong} refers to the fact that for any congruence on a finite lattice, each congruence class has a unique minimal element.
(See Section~\ref{weak sec}).
The corresponding quotient lattice is isomorphic to the subposet induced by the set of elements which are minimal in their congruence class. 
In the course of proving Theorem~\ref{cong}, we also establish the following results.  (The map $w\mapsto ww_0$ in Proposition~\ref{anti} is explained in Section~\ref{sort sec}.)

\begin{theorem}
\label{sublattice}
Let~$c$ be a Coxeter element of a finite Coxeter group~$W\!$.\@
The \mbox{$c$-sortable} elements constitute a sublattice of the weak order on $W\!$.\@
\end{theorem}

\begin{prop}
\label{anti}
The map $w\mapsto ww_0$ maps the congruence $\Theta_c$ to the congruence $\Theta_{c^{-1}}$.
In particular, the weak order on \mbox{$c$-sortable} elements is anti-isomorphic to the weak order on $c^{-1}$-sortable elements.
\end{prop}

The second main result of the paper concerns the connection between sortable elements and the Cambrian lattices of~\cite{cambrian}. 
These lattices were originally proposed as a simple construction valid for arbitrary Coxeter groups which, in the cases of types A and B, was known to reproduce the combinatorics and geometry of the generalized associahedra.
However, the lattice-theoretic definition does not immediately shed light on the combinatorics of the lattice, so that in particular the definition was not proven, outside of types A and B and small examples, to relate to the combinatorics of \mbox{$W$-Catalan} numbers.
The content of the following theorem is that sortable elements provide a concrete combinatorial realization of the Cambrian lattices for any $W\!$.\@

\begin{theorem}
\label{camb}
The congruence $\Theta_c$ is the Cambrian congruence associated to~$c$.  
In particular, the Cambrian lattice associated to~$c$ is the weak order on \mbox{$c$-sortable} elements.
\end{theorem}
The proof of Theorem~\ref{camb} is accomplished using the geometric model for lattice congruences of the weak order laid out in~\cite{congruence}.
The results of this paper prove some of the conjectures of~\cite{cambrian}.
Many of the remaining conjectures are proved in~\cite{camb_fan}.
In light of Theorem~\ref{camb}, Proposition~\ref{anti} is equivalent to \cite[Theorem~3.5]{cambrian}.

We conclude the introduction with two examples which illustrate the results of the paper.
The definitions underlying these examples appear in later sections.

\begin{figure}[ht]
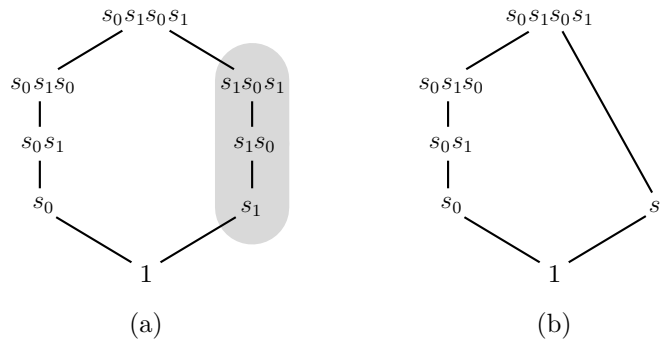

\centerline{
\begin{tabular}{cccccc}
\scalebox{.8}{\epsfbox{B2cong.ps}}&&&&&\scalebox{.8}{\epsfbox{B2camb.ps}}\\[8 pt]
(a)&&&&&(b)
\end{tabular}
}
\caption{An example of $\Theta_c$ and the associated Cambrian lattice}
\label{B2camb}
\end{figure}

\begin{example}\rm
\label{B2lat}
Consider the case $W=B_2$ with $S=\set{s_0,s_1}$, $m(s_0,s_1)=4$ and $c=s_0s_1$.
Figure~\ref{B2camb}.a shows the congruence $\Theta_c$ on the weak order on~$W\!$.\@  
The shaded $3$-element chain is a congruence class and each other congruence class is a singleton.
Figure~\ref{B2camb}.b shows the subposet of the weak order induced by the sortable elements, or equivalently the bottom elements of congruence classes.
The map $w\mapsto ww_0$ acts on the Hasse diagram in Figure~\ref{B2camb}.a by rotating through a half-turn.
One easily verifies Theorem~\ref{sublattice} and Proposition~\ref{anti} in this example.
\end{example}

\begin{example}\rm
\label{A3lat}
A more substantive example is provided by $W=A_3$ with $S=\set{s_1,s_2,s_3}$, $m(s_1,s_2)=m(s_2,s_3)=3$, $m(s_1,s_3)=2$ and $c=s_2s_1s_3$.
This Coxeter group is isomorphic to the symmetric group $S_4$ as generated by the simple transpositions $s_i=(i\,\,i\!+\!1)$.
Figure~\ref{sortcong} shows the congruence $\Theta_c$ on the weak order on~$A_3$.
Each element is represented by its $c$-sorting word (see Section~\ref{sort sec}) including the inert dividers ``$|$'', with the symbols $s_i$ replaced by $i$ throughout.
(The identity element is represented by the empty word.)
Non-trivial $\Theta_c$-classes are indicated by shading, and each unshaded element is the unique element in its congruence class.
The antiautomorphism $w\mapsto ww_0$ corresponds to rotating the diagram through a half-turn.

\begin{figure}[ht]

\begin{tabular}{c}
\scalebox{0.7}{\epsfbox{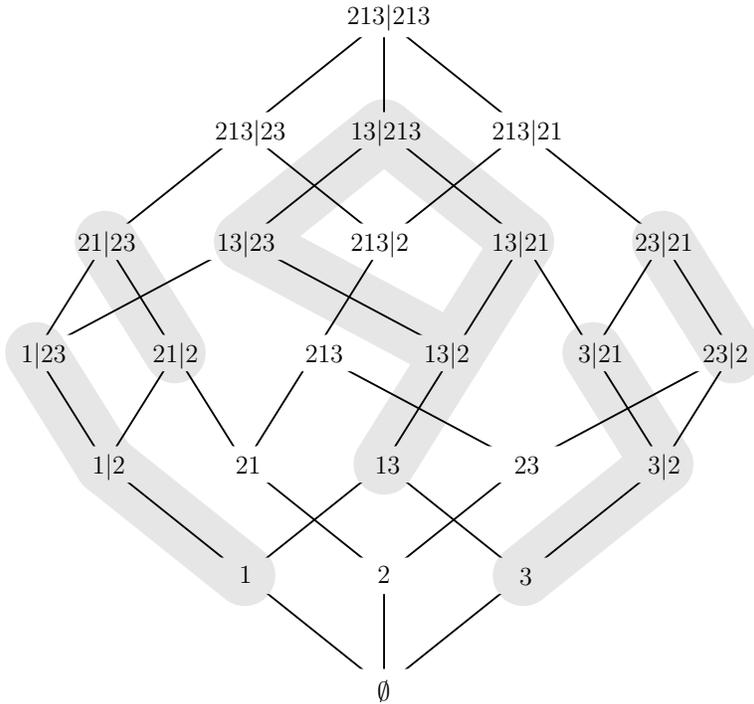}}\\[-7pt]
$\emptyset$
\end{tabular}

\caption{Another example of $\Theta_c$}
\label{sortcong}
\end{figure}

The content of Theorem~\ref{camb} in this example is that $\Theta_c$ is the smallest congruence of the weak order on~$A_3$ which sets $s_1\equiv s_1s_2$ and $s_3\equiv s_3s_2$.
Combining Theorem~\ref{camb} with Proposition~\ref{anti}, we obtain the assertion that $\Theta_c$ is the smallest congruence which sets $s_1s_3s_2s_3\equiv s_1s_3s_2s_1s_3\equiv s_1s_3s_2s_1$.
\end{example}

\hfill
\newpage

\section{Sortable elements}
\label{sort sec}
In this section we quickly review the definition and first properties of Coxeter groups.
(More detail, including proofs of assertions not proven here, can be found in one of the standard references \cite{Bj-Br,Bourbaki,Humphreys}.)
We then review the definition of sortable elements and quote and prove results which are used in later sections.

A {\em Coxeter group}~$W$ is a group with a presentation as the group generated by~a set~$S$ subject to the relations $(st)^{m(s,t)}=1$ for $s,t\in S$.
Here $m(s,s)=1$ and $m(s,t)\ge 2$ for $s\neq t$, with $m(s,t)=\infty$ meaning that no relation of the form $(st)^m$ holds.
It can be shown that $m(s,t)$ is the order of $st$ in~$W$.
(\mbox{\em A priori}, we only know that the order of $st$ divides $m(s,t)$.)
The {\em rank} of~$W$ is~$|S|$.
A group can have more than one non-equivalent presentation of this form but, as usual, we take the term ``Coxeter group'' to imply a distinguished choice of~$S$, called the {\em simple generators} of~$W$.
The presentation of~$W$ is encoded in the {\em Coxeter diagram} of~$W$.
This is a graph whose vertices are the simple generators, with edges $s\,$---$\,t$ whenever $m(s,t)>2$.
Edges are labeled by the number $m(s,t)$ except that, by convention, if $m(s,t)=3$ then the edge is left unlabeled.

Each $w\in W$ can be written, in many different ways, as a word in the alphabet~$S$.
The smallest length of a word for~$w$ is called the {\em length} of~$w$, denoted $l(w)$.
A word of length $l(w)$ representing~$w$ is called {\em reduced}.
For any $s\in S$, $l(sw)=l(w)\pm 1$ and $l(sw)<l(w)$ if and only if there is some reduced word for~$w$ starting with the letter~$s$.
The analogous statement holds for $l(ws)$ and $l(w)$.

Every Coxeter group has a {\em reflection representation}:  a representation as a group generated by orthogonal reflections of a real vector space.
A Coxeter group is finite if and only if it has such a representation in a vector space with a Euclidean inner product.
An element $t\in W$ acts as an orthogonal reflection in such a representation if and only if it is conjugate to a simple generator.
Thus the set $T=\set{wsw^{-1}:w\in W, s\in S}$ is called the set of {\em reflections} in~$W$.
The {\em (left) inversion set} of $w\in W$ is $\set{t\in T:l(tw)<l(w)}$.  
The length $l(w)$ is equal to $|I(w)|$.
An element of~$W$ is uniquely determined by its inversion set.

The {\em (right) weak order} is the partial order on~$W$ whose cover relations are $w\covered ws$ whenever $l(w)<l(ws)$.
Equivalently, $v\le w$ if and only if $I(v)\subseteq I(w)$.
Throughout this paper, ``$\le$'' denotes the weak order and ``$W$'' denotes both the group~$W$ and the partially ordered set~$W$.
This partial order is a meet-semilattice in general, and a lattice exactly when~$W$ is finite~\cite{orderings}.
For a simple generator~$s$ and an interval $[u,v]$ in the weak order, if $l(su)<l(u)$ then the involution $w\mapsto sw$ is an isomorphism between the intervals $[su,sv]$ and $[u,v]$.
The minimal element of~$W$ is the identity $1$ and when~$W$ is finite, the longest element $w_0$ is the unique maximal element of~$W$.
In this case, the map $w\mapsto ww_0$ is an antiautomorphism of weak order, because $I(ww_0)=T-I(w)$.

A {\em Coxeter element}~$c$ of~$W$ is an element represented by a (necessarily reduced) word $a_1a_2\cdots a_n$ where $S=\set{a_1,\ldots,a_n}$ and $n=|S|$.
Fix~$c$ and a particular word $a_1a_2\cdots a_n$ for~$c$ and write a half-infinite word 
\[c^\infty=a_1a_2\cdots a_n|a_1a_2\cdots a_n|a_1a_2\cdots a_n|\ldots\]
The symbols ``$|$'' are inert ``dividers'' which facilitate the definition of sortable elements.
When subwords of $c^\infty$ are interpreted as expressions for elements of $W,$ the dividers are ignored.
The {\em $c$-sorting word} for $w\in W$ is the lexicographically first (as a sequence of positions in $c^\infty$) subword of $c^\infty$ which is a reduced word for~$w$.
The $c$-sorting word can be interpreted as a sequence of subsets of~$S$:
Each subset in the sequence is the set of letters of the $c$-sorting word which occur between two adjacent dividers.

An element $w\in W$ is {\em \mbox{$c$-sortable}} if its $c$-sorting word defines a sequence of subsets which is weakly decreasing under inclusion.
Since any two reduced words for~$c$ are related by commutation of letters, the $c$-sorting words for~$w$ arising from different reduced words for~$c$ are related by commutations of letters, with no commutations across dividers.
In particular, the set of \mbox{$c$-sortable} elements does not depend on the choice of reduced word for~$c$.
Examples~\ref{B2lat} and~\ref{A3lat} in Section~\ref{intro} illustrate the definition of sortable elements.
(In Example~\ref{A3lat}, the \mbox{$c$-sortable} elements are the elements at the bottom of their congruence class, including elements which are unique in their class.)

For any $J\subseteq S$, let $W_J$ be the subgroup of~$W$ generated by $J$.
Such subgroups are called {\em standard parabolic subgroups}.
For any $w\in W,$ there is a unique factorization $w=w_J\cdot\closeleftq{J}{w}$ such that $w_J\in W_J$ and $\closeleftq{J}{w}$ has $l(s\cdot\closeleftq{J}{w})>l(\closeleftq{J}{w})$ for all $s\in J$.
The element $w_J$ has inversion set $I(w_J)=I(w)\cap W_J$.
The element $w_J\cdot\closeleftq{J}{(w_0)}$ has inversion set $I(w)\cup\set{t\in T:t\not\in W_J}$.
The map $w\mapsto w_J$ is a lattice homomorphism \cite[Section~6]{congruence} for any $J\subseteq S$.
(That is, $(x\join y)_J=(x_J)\join(y_J)$ and similarly for meets.)
This map is also compatible with the antiautomorphism $w\mapsto ww_0$ in the sense that $(ww_0)_J=w_J(w_0)_J$.
The set $W_J$ is a lower interval in the poset~$W$ with maximal element~$(w_0)_J$.
Most often, the subset $J$ is $\br{s}:=S-\set{s}$ for some element~$s$ of~$S$.

An {\em initial letter} of a Coxeter element~$c$ is a simple generator which is the first letter of some reduced word for~$c$.
Similarly a {\em final letter} of~$c$ is a simple generator which occurs as the last letter of some reduced word for~$c$.
If~$s$ is initial in~$c$ then $scs$ is a Coxeter element and~$s$ is final in $scs$.
The next two lemmas \cite[Lemmas 2,4 and 2.5]{sortable} constitute an inductive characterization of sortable elements.
As a base for the induction, $1$ is \mbox{$c$-sortable} for any~$c$.

\begin{lemma}
\label{sc}
Let~$s$ be an initial letter of~$c$ and let $w\in W$ have $l(sw)>l(w)$.
Then~$w$ is \mbox{$c$-sortable} if and only if it is an $sc$-sortable element of~$W_{\br{s}}$.
\end{lemma}

\begin{lemma}
\label{scs}
Let~$s$ be an initial letter of~$c$ and let $w\in W$ have $l(sw)<l(w)$.
Then~$w$ is \mbox{$c$-sortable} if and only if $sw$ is $scs$-sortable.
\end{lemma}

For the rest of the paper we confine our attention to the case where~$W$ is a finite Coxeter group.
In particular, all proofs should be assumed to apply only to the finite case.

The remainder of the section is devoted to quoting or proving preliminary results.
For any $J\subseteq S$, the {\em restriction} of a Coxeter element~$c$ to~$W_J$ is the Coxeter element for~$W_J$ obtained
by deleting the letters $S-J$ from any reduced word for~$c$.
The next lemma is immediate from the definition of sortable elements and the proposition following it is \cite[Corollary~4.5]{sortable}.

\begin{lemma}
\label{sort para easy}
Let~$c$ be a Coxeter element of~$W,$ let $J\subseteq S$ and let~$c'$ be the Coxeter element of~$W_J$ obtained by restriction.
If $w\in W_J$ is $c'$-sortable then~$w$ is \mbox{$c$-sortable} as an element of $W\!$.\@
\end{lemma}

\begin{prop}
\label{sort para}
Let~$c$ be a Coxeter element of~$W,$ let $J\subseteq S$ and let~$c'$ be the Coxeter element of~$W_J$ obtained by restriction.
If~$w$ is \mbox{$c$-sortable} then~$w_J$ is $c'$-sortable.
\end{prop}

A {\em cover reflection} of $w\in W$ is a reflection $t\in T$ such that $tw=ws$ with $l(ws)<l(w)$.
The term ``cover reflection'' refers to the fact that $w\covers ws$ is a cover relation in the weak order.
Equivalently, $t\in T$ is a cover reflection of~$w$ if and only if $I(w)-\set{t}$ is the inversion set of some $w'\in W$.
In this case $w'=tw\covered w$.
Let $\cov(w)$ denote the set of cover reflections of~$w$.
The following proposition is one direction of the rephrasing of \cite[Theorem~6.1]{sortable} described in \cite[Remark 6.9]{sortable}.

\begin{prop}
\label{nc cov}
If~$x$ and~$y$ are \mbox{$c$-sortable} and $\cov(x)=\cov(y)$ then $x=y$.
\end{prop}

The next lemma follows from \cite[Lemma~5.2]{sortable} and \cite[Theorem~6.1]{sortable} as described in \cite[Remark 6.9]{sortable}.

\begin{lemma}
\label{s cov}
If~$s$ is initial in~$c$ and $x\in W_{\br{s}}$ is $sc$-sortable then there exists a \mbox{$c$-sortable} element~$w$ with $\cov(w)=\set{s}\cup\cov(x)$.
\end{lemma}

We conclude with four lemmas on the join operation as it relates to cover reflections and sortable elements.
The first two lemmas are used only in the proof of the last two, which are applied in Section~\ref{weak sec}.

\begin{lemma}
\label{s join w br s}
For $w\in W,$ if~$s$ is a cover reflection of~$w$ and every other cover reflection of~$w$ is in~$W_{\br{s}}$ then $w=s\join w_{\br{s}}$.
\end{lemma}
\begin{proof}
Since~$s$ is a cover reflection of~$w$, it is in particular an inversion of~$w$, and since $I(s)=\set{s}$, we have $s\le w$.
Any element~$w$ has $w_{\br{s}}\le w$, so~$w$ is an upper bound for $s$ and $w_{\br{s}}$.
Since~$s$ is a cover reflection of~$w$ and every other cover reflection is in~$W_{\br{s}}$, any element~$x$ covered by~$w$ has either $s\not\le x$ or $w_{\br{s}}\not\le x$ (the latter because $I(w_{\br{s}})=I(w)\cap W_{\br{s}}$).
Thus $w=s\join w_{\br{s}}$.
\end{proof}

\begin{lemma}
\label{cov w br s}
For $x\in W_{\br{s}}$, $\cov(s\join x)=\cov(x)\cup\set{s}$.
\end{lemma}
\begin{proof}
Let $w=s\join x$.
First, we show that the reflection~$s$ is a cover reflection of~$w$.
If not, then let $w'$ be any element covered by~$w$ and weakly above~$x$.
Since~$s$ is not a cover reflection of~$w$ we have $s\in I(w')$, so $w'$ is above both~$s$ and~$x$, contradicting the fact that $w=s\join x$.

Now let $t\neq s$ be a cover reflection of~$w$, so that $I(w)-I(tw)=\set{t}$.
If $t\not\in W_{\br{s}}$ then since $I(x)\subset I(w)$ and $I(x)\subseteq W_{\br{s}}$ we also have $I(x)\subseteq I(tw)$.
Since $s\neq t$ we have $s\le tw$, and this contradicts the fact that $w=s\join x$.
This contradiction proves that $t\in W_{\br{s}}$.
Now, $I(w_{\br{s}})=I(w)\cap W_{\br{s}}$ and $I((tw)_{\br{s}})=I(w)\cap W_{\br{s}}$ and since $t\in W_{\br{s}}$, $I(w_{\br{s}})-I((tw)_{\br{s}})=\set{t}$.
Thus $t$ is a cover reflection of $w_{\br{s}}$.
Applying the lattice homomorphism $y\mapsto y_{\br{s}}$ to the equality $w=s\join x$, we have $w_{\br{s}}=x$.

Conversely suppose~$x$ covers $tx$.
Let~$y$ be any element covered by~$w$ and weakly above $s\join tx$.
Then $y=t'w$ for some reflection $t'$.
Since $I(x)\subset I(w)$, $I(x)\not\subset I(y)$ and $I(w)-I(y)=\set{t'}$, we must have $t'\in I(x)$.
Since $I(tx)\subset I(s\join tx)\subset I(y)$ we must have $t'\not\in I(tx)$.
But $I(x)-I(tx)=\set{t}$ so $t=t'$.
Thus $t$ is a cover reflection of~$w$.
\end{proof}

\begin{lemma}
\label{s join 1}
If~$s$ is initial in~$c$ and $x\in W_{\br{s}}$ is $sc$-sortable then $s\join x$ is \mbox{$c$-sortable}.
\end{lemma}
\begin{proof}
Let~$s$ be initial in~$c$, let $x\in W_{\br{s}}$ be $sc$-sortable and let~$w$ be the \mbox{$c$-sortable} element, as in Lemma~\ref{s cov}, such that $\cov(w)=\set{s}\cup\cov(x)$.
By Lemma~\ref{s join w br s}, $w=s\join w_{\br{s}}$.
By Proposition~\ref{sort para}, $w_{\br{s}}$ is $sc$-sortable and by Lemma~\ref{cov w br s}, $\cov(w_{\br{s}})=\cov(w)-\set{s}=\cov(x)$.
But Proposition~\ref{nc cov} says that $w_{\br{s}}=x$ and thus in particular $s\join x=w$ is \mbox{$c$-sortable}.
\end{proof}

\begin{lemma}
\label{s join 2}
If~$s$ is final in~$c$ and $w\in W$ is \mbox{$c$-sortable} with $l(sw)<l(w)$ then $w=w_{\br{s}}\join s$.
\end{lemma}
\begin{proof}
By \cite[Lemmas 6.6 and 6.7]{sortable}, $s$ is a cover reflection of~$w$ and every other cover reflection is in $W_{\br{s}}$.
Thus Lemma~\ref{s join w br s} says that $w=w_{\br{s}}\join s$.
\end{proof}

\section{Sortable elements and lattice congruences}
\label{weak sec}
In this section we review the definition and an order-theoretic characterization of lattice congruences and prove Theorem~\ref{cong}, Theorem~\ref{sublattice} and Proposition~\ref{anti}.
We remind the reader that all proofs given here are valid only in the finite case.

A {\em congruence} on a lattice $L$ is an equivalence relation $\Theta$ on $L$ such that whenever $a_1\equiv b_1$ and $a_2\equiv b_2$ then $a_1\meet a_2\equiv b_1\meet b_2$ and $a_1\join a_2\equiv b_1\join b_2$.  
The quotient of $L$ mod $\Theta$ is a lattice defined on the $\Theta$-congruence classes.
Denote the $\Theta$-congruence class of $a$ by $[a]_\Theta$ and set $[a_1]_\Theta\meet[a_2]_\Theta=[a_1\meet a_2]_\Theta$ and $[a_1]_\Theta\join[a_2]_\Theta=[a_1\join a_2]_\Theta$. 

The following order-theoretic characterization of congruences of a finite lattice was introduced in~\cite{dissective}.
It is a straightforward exercise to prove the characterization, which also follows from a characterization of congruences of general lattices due to Chajda and Sn\'a\v{s}el~\cite{Cha-Sn}.

\pagebreak[3]

\begin{prop}
\label{order cong}
When $L$ is finite, an equivalence relation on $L$ is a lattice congruence if and only if it has the following three properties.
\begin{enumerate}
\item[(i) ] Every equivalence class is an interval.
\item[(ii) ] The downward projection $\pidown:L\rightarrow L$, mapping each element to the minimal element in its equivalence class, is order-preserving.
\item[(iii) ] The upward projection $\piup:L\rightarrow L$, mapping each element to the maximal element in its equivalence class, is order-preserving.
\end{enumerate}
Furthermore the quotient of $L$ mod $\Theta$ is isomorphic to the subposet of $L$ induced by the set $\pidown L$ of bottom elements of congruence classes and $\pidown$ is a homomorphism from $L$ to $\pidown L$.
\end{prop}

The proof of Theorem~\ref{cong} uses the characterization of congruences in Proposition~\ref{order cong}.
Specifically, we construct order-preserving maps $\pidown^c$ and $\piup_c$ such that
first, $\pidown^c(x)=x$ if and only if~$x$ is \mbox{$c$-sortable} and
second, the equivalence setting $x\equiv y$ if and only if $\pidown^c(x)=\pidown^c(y)$ has equivalence classes of the form $[\pidown^c(x),\piup_c(x)]$.

Some lemmas which contribute to the proof make use of induction on length and rank in a way that is similar to the proofs presented in~\cite{sortable}.
In the course of these inductive proofs, we repeatedly apply the following fact: 
For $s\in S$ and $x,y\in W,$ if $l(sx)<l(x)$ and $l(sy)<l(y)$ then $x\le y$ if and only if $sx\le sy$.

We also make use of the following notational convention:
The explicit reference to~$c$ in $\pidown^c$ and $\piup_c$ is used as a way of specifying a standard parabolic subgroup.
For example, if~$s$ is initial in~$c$ then $sc$ is a Coxeter element in the standard parabolic subgroup~$W_{\br{s}}$.
Thus the notation $\pidown^{sc}$ refers to a map on~$W_{\br{s}}$.

We now construct a projection $\pidown^c$.
This is done by induction on the rank of~$W$ and on the length of the element to which $\pidown^c$ is applied.
Define $\pidown^c(1)=1$ and for any initial letter~$s$ of~$c$, define
\begin{equation}
\label{pidown def}
\pidown^c(w)=\left\lbrace\begin{array}{ll}
s\cdot\pidown^{scs}(sw)&\mbox{if }l(sw)<l(w),\mbox{ or}\\
\pidown^{sc}(w_{\br{s}})&\mbox{if }l(sw)>l(w).
\end{array}\right.
\end{equation}
As written, each step of this inductive definition depends on a choice of an initial letter of a Coxeter element.
The following proposition implies that $\pidown^c$ is well-defined and that $\pidown^c\circ\pidown^c=\pidown^c$.

\begin{prop}
\label{pidown max}
For any $w\in W,$ $\pidown^c(w)$ is the unique maximal \mbox{$c$-sortable} element weakly below~$w$.
\end{prop}
\begin{proof}
Let $w\in W$ and let~$s$ be initial in~$c$.
If $l(sw)<l(w)$ then by induction $\pidown^{scs}(sw)$ is the unique maximal $scs$-sortable element weakly below $sw$. 
Then \mbox{$s\cdot\pidown^{scs}(sw)$} is weakly below~$w$, and by Lemma~\ref{scs}, $s\cdot\pidown^{scs}(sw)$ is \mbox{$c$-sortable}.
Let~$x$ be any \mbox{$c$-sortable} element below~$w$.  
If $l(sx)<l(x)$ then $sx$ is an $scs$-sortable element below $sw$, and therefore $sx$ is below $\pidown^{scs}(sw)$, so that~$x$ is below
$s\cdot\pidown^{scs}(sw)=\pidown^c(w)$.
If $l(sx)>l(x)$ then by Lemmas~\ref{sc} and~\ref{s join 1}, $x\join s$ is \mbox{$c$-sortable}.
Now $s\le w$ so $x\join s\le w$ and since $x\join s$ is shortened on the left by~$s$, by the previous case $x\join s$ is below $\pidown^c(w)$, and therefore~$x$ is below $\pidown^c(w)$.

If $l(sw)>l(w)$ then by Lemma~\ref{sc}, any \mbox{$c$-sortable} element~$x$ below~$w$ is an $sc$-sortable element of~$W_{\br{s}}$.
In particular, $x\le w_{\br{s}}$.
By induction on the rank of~$W,$ $\pidown^c(w)=\pidown^{sc}(w_{\br{s}})$ is the unique maximal such.
\end{proof}

\begin{cor}
\label{pidown op}
The map $\pidown^c$ is order-preserving on the weak order.
\end{cor}
\begin{proof}
Suppose $x\le y$.
By Proposition~\ref{pidown max}, $\pidown^c(x)$ is a \mbox{$c$-sortable} element below~$x$ and therefore below~$y$.
Thus Proposition~\ref{pidown max} says that $\pidown^c(x)\le\pidown^c(y)$.
\end{proof}

We now prove Theorem~\ref{sublattice}, which asserts that \mbox{$c$-sortable} elements constitute a sublattice of the weak order.

\begin{proof}[Proof of Theorem~\ref{sublattice}]
Let~$x$ and~$y$ be \mbox{$c$-sortable}.
We first show that $x\meet y$ is also \mbox{$c$-sortable}.
Choose an initial letter~$s$ of~$c$.
If $l(sx)>l(x)$ and $l(sy)>l(y)$ then~$x$ and~$y$ are both in~$W_{\br{s}}$.
By induction on the rank of~$W,$ $x\meet y$ is an $sc$-sortable element of~$W_{\br{s}}$, so it is also \mbox{$c$-sortable} by Lemma~\ref{sort para easy}.

If exactly one of~$x$ and~$y$ are shortened on the left by~$s$, we may as well take $l(sx)<l(x)$ and $l(sy)>l(y)$.
Then $y\in W_{\br{s}}$.
Since $W_{\br{s}}$ is a lower interval in $W$ and $w\mapsto w_{\br{s}}$ is a lattice homomorphism
\[x\meet y=(x\meet y)_{\br{s}}=x_{\br{s}}\meet y_{\br{s}}=x_{\br{s}}\meet y.\]
By Proposition~\ref{sort para}, $x_{\br{s}}$ is $sc$-sortable, so by the previous case $x\meet y$ is \mbox{$c$-sortable}.

If $l(sx)<l(x)$ and $l(sy)<l(y)$ then by Lemma~\ref{scs}, $sx$ and $sy$ are both $scs$-sortable.
By induction on length, $sx\meet sy$ is $scs$-sortable.
Since left multiplication by~$s$ is an isomorphism from the interval $[s,w_0]$ to the interval $[1,sw_0]$, $sx\meet sy$ is lengthened on the 
left by~$s$ and $x\meet y=s(sx\meet sy)$.
Now Lemma~\ref{scs} says that $s(sx\meet sy)$ is \mbox{$c$-sortable}.

We now show that $x\join y$ is \mbox{$c$-sortable}.
Since $x\join y\ge x$ we have $\pidown^c(x\join y)\ge\pidown^c(x)$ by Corollary~\ref{pidown op}.
By Proposition~\ref{pidown max} $\pidown^c(x)=x$.
Similarly $\pidown^c(x\join y)\ge y$, so $\pidown^c(x\join y)$ is an upper bound for~$x$ and~$y$.
By definition of join, $x\join y\le \pidown^c(x\join y)$, but Proposition~\ref{pidown max} says that $\pidown^c(x\join y)\le x\join y$.
Thus $x\join y=\pidown^c(x\join y)$, which is \mbox{$c$-sortable}.
\end{proof}

We continue towards a proof of Theorem~\ref{cong} by defining the upward projection corresponding to $\pidown^c$.

Call $w\in W$ {\em $c$-antisortable} if~$ww_0$ is $c^{-1}$-sortable.
Define an upward projection map $\piup_c$ by setting $\piup_c(w)=\left(\pidown^{(c^{-1})}(ww_0) \right)w_0$.
Since $w\mapsto ww_0$ is an antiautomorphism of the right weak order, it is immediate from Proposition~\ref{pidown max} and Corollary~\ref{pidown op} that $\piup_c(w)$ is the unique minimal element among $c$-antisortable elements above~$w$, that $\piup_c\circ\piup_c=\piup_c$ and that $\piup_c$ is order-preserving.

\begin{lemma}
\label{piup formula}
If~$s$ is a final letter of~$c$ then
\begin{equation*}
\piup_c(w)=\left\lbrace\begin{array}{ll}
s\cdot\piup_{scs}(sw)&\mbox{if }l(sw)>l(w),\mbox{ or}\\
\piup_{cs}(w_{\br{s}})\cdot\leftq{\br{s}}{w_0}&\mbox{if }l(sw)<l(w).
\end{array}\right.
\end{equation*}
\end{lemma}
\begin{proof}
If $l(sw)>l(w)$ then $l(sww_0)<l(ww_0)$, so 
\[\piup_c(w)=s\cdot\left(\pidown^{(sc^{-1}s)}(sww_0)\right)w_0=s\cdot\piup_{scs}(sw).\]
If $l(sw)<l(w)$ then $l(sww_0)>l(ww_0)$, so $\piup_c(w)=\pidown^{(sc^{-1})}\left((ww_0)_{\br{s}}\right)w_0$.
But $(ww_0)_{\br{s}}=w_{\br{s}}(w_0)_{\br{s}}$, so
\[\pidown^{(sc^{-1})}\left((ww_0)_{\br{s}}\right)w_0=\pidown^{(sc^{-1})}\left(w_{\br{s}}(w_0)_{\br{s}}\right)(w_0)_{\br{s}}\cdot\leftq{\br{s}}{w_0}=\piup_{cs}(w_{\br{s}})\cdot\leftq{\br{s}}{w_0}.\]
\end{proof}

In order to construct the congruence $\Theta_c$, it is necessary to relate the fibers of $\pidown^c$ to those of $\piup_c$.
This is done by induction on rank and length as in the proofs of Proposition~\ref{pidown max} and Theorem~\ref{sublattice}.
However, the recursive definition of $\pidown^c$ requires an initial letter of $c$ and the corresponding property of $\piup_c$ (Lemma~\ref{piup formula}) requires a final letter of $c$.
Thus an additional tool is needed, and this tool is provided by the dual (Lemma~\ref{piup alt form}) of the following lemma.

\begin{lemma}
\label{pidown alt form}
If~$s$ is a final letter of~$c$ and $l(sw)<l(w)$ then $\pidown^c(w)=s\join\pidown^{cs}(w_{\br{s}})$.
\end{lemma}
\begin{proof}
Let $x=s\join\pidown^{cs}(w_{\br{s}})$.
By Lemma~\ref{sort para easy} and Proposition~\ref{pidown max}, $\pidown^{cs}(w_{\br{s}})$ is a \mbox{$c$-sortable} element.
Thus both~$s$ and $\pidown^{cs}(w_{\br{s}})$ are \mbox{$c$-sortable} elements below~$w$ and therefore~$x$ is a \mbox{$c$-sortable} element below~$w$ by Theorem~\ref{sublattice}.
(Note that the weaker result, Lemma~\ref{s join 1}, applies to an initial letter, and thus cannot be used to show that~$x$ is \mbox{$c$-sortable}.)
By Lemma~\ref{pidown max}, $x\le\pidown^c(w)$.
Since $l(sx)<l(x)$, $\pidown^c(w)$ is also shortened on the left by~$s$.
Now Lemma~\ref{s join 2} says that $\pidown^c(w)=s\join\left(\pidown^c(w)_{\br{s}}\right)$.
Because $\pidown^c(w)\le w$, we have $\pidown^c(w)_{\br{s}}\le w_{\br{s}}$.
Since $\pidown^c(w)_{\br{s}}$ is a $cs$-sortable element below~$w_{\br{s}}$, it is below $\pidown^{cs}(w_{\br{s}})$ by Proposition~\ref{pidown max}.
Therefore $x=s\join\pidown^{cs}(w_{\br{s}})\ge s\join\left(\pidown^c(w)_{\br{s}}\right)=\pidown^c(w)$.
\end{proof}

\begin{lemma}
\label{piup alt form}
If~$s$ is an initial letter of~$c$ and $l(sw)>l(w)$ then $\piup_c(w)=sw_0\meet\left(\piup_{sc}(w_{\br{s}})\cdot\leftq{\br{s}}{w_0}\right)$.
\end{lemma}
\begin{proof}
Starting with Lemma~\ref{pidown alt form}, replace~$w$ by $ww_0$ and $c$ by $c^{-1}$, multiply both sides on the right by $w_0$ and apply the fact that $w\mapsto ww_0$ is an antiautomorphism.
\end{proof}

We now relate $\pidown^c$ to $\piup_c$.

\begin{prop}
\label{equiv}
The maps $\pidown^c$ and $\piup_c$ are compatible in the following senses.
\begin{enumerate}
\item[(i) ]For any $x,y\in W,$ $\pidown^c(x)=\pidown^c(y)$ if and only if $\piup_c(x)=\piup_c(y)$.
\item[(ii) ]$\piup_c\circ\pidown^c=\piup_c$ and $\pidown^c\circ\piup_c=\pidown^c$.
\end{enumerate}
\end{prop}
\begin{proof}
By the antisymmetry $w\mapsto ww_0$, to prove (i) it suffices to prove the ``only if'' direction.
We treat first the special case of (i) where $x\le y$.
Suppose $\pidown^c(x)=\pidown^c(y)$ and $x\le y$.
Let~$s$ be initial in~$c$. 
If $l(sx)<l(x)$ then since $x\le y$, $l(sy)<l(y)$.
Thus $\pidown^c(x)=s\cdot\pidown^{scs}(sx)$ and $\pidown^c(y)=s\cdot\pidown^{scs}(sy)$, so $\pidown^{scs}(sx)=\pidown^{scs}(sy)$.
Since $sx\le sy$, by induction on $l(x)$, $\piup_{scs}(sx)=\piup_{scs}(sy)$.
By Lemma~\ref{piup formula} (with $c$ replaced by $scs$ and~$w$ replaced by~$x$ or~$y$) $\piup_{scs}(sx)=s\cdot\piup_c(x)$ and $\piup_{scs}(sy)=s\cdot\piup_c(y)$ so $\piup_c(x)=\piup_c(y)$.

If $l(sx)>l(x)$, we claim that $l(sy)>l(y)$.
If not then $\pidown^c(x)=\pidown^{sc}(x_{\br{s}})$ and $\pidown^c(y)=s\cdot\pidown^{scs}(sy)$.
In particular, $\pidown^c(x)$ is lengthened on the left by~$s$ but $\pidown^c(y)$ is shortened on the left by~$s$, contradicting the supposition that $\pidown^c(x)=\pidown^c(y)$.
This contradiction proves the claim.
Thus $\pidown^c(x)=\pidown^{sc}(x_{\br{s}})$ and $\pidown^c(y)=\pidown^{sc}(y_{\br{s}})$, so that $\pidown^{sc}(x_{\br{s}})=\pidown^{sc}(y_{\br{s}})$.
Since $x_{\br{s}}\le y_{\br{s}}$, by induction on the rank of~$W$ $\piup_{sc}(x_{\br{s}})=\piup_{sc}(y_{\br{s}})$.
By Lemma~\ref{piup alt form}, $\piup_c(x)=\piup_c(y)$.

Having established (i) in the case $x\le y$, we now prove (ii).
Any $y\in W$ has $\pidown^c(y)\le y$ and $\pidown^c(\pidown^c(y))=\pidown^c(y)$, so setting $x=\pidown^c(y)$ in the special case of (i) already proved, $\piup_c(\pidown^c(y))=\piup_c(y)$.
Thus $\piup_c\circ\pidown^c=\piup_c$.
The antisymmetry $w\mapsto ww_0$ now implies that $\pidown^c\circ\piup_c=\pidown^c$ as well.

Finally, we prove the general case of (i).
If~$x$ and~$y$ are unrelated in the weak order and $\pidown^c(x)=\pidown^c(y)$ then 
$\piup_c(x)=\piup_c(\pidown^c(x))=\piup_c(\pidown^c(y))=\piup_c(y)$.
\end{proof}

We now prove the main theorem of the section.

\begin{proof}[Proof of Theorem~\ref{cong}]
For each $w\in W,$ let $D(w)=\set{y\in W:\pidown^c(y)=\pidown^c(w)}$ and let $U(w)=\set{y\in W:\piup_c(y)=\piup_c(w)}$.
Since $\pidown^c$ is order-preserving, $\pidown^c(w)$ is the unique minimal element of $D(w)$, and similarly, $\piup_c(w)$ is the unique maximal element of $U(w)$.
Proposition~\ref{equiv} states that $D(w)=U(w)$, and since $\pidown^c$ is order preserving, $D(w)$ is the entire interval $[\pidown^c(w),\piup_c(w)]$ in the right weak order.

Thus the fibers of the map $\pidown^c$ form a decomposition of the weak order on~$W$ into intervals and $\pidown^c$ and  $\piup_c$ are the order-preserving maps required in Proposition~\ref{order cong}.
Therefore these intervals are the congruence classes of some congruence $\Theta_c$.
\end{proof}

The preceding considerations also prove Proposition~\ref{anti}, as we now explain.
The congruence classes of $\Theta_c$ are of the form $[\pidown^c(w),\piup_c(w)]$.
Applying the anti-automorphism $w\mapsto ww_0$ we obtain intervals of the form $[\pidown^{(c^{-1})}(ww_0),\piup_{(c^{-1})}(ww_0)]$.
Thus the anti-automorphism $w\mapsto ww_0$ takes $\Theta_c$ to $\Theta_{c^{-1}}$.

\begin{remark}\rm
The \mbox{$c$-sortable} elements of a finite Coxeter group~$W$ are counted by the \mbox{$W$-Catalan} number, as is proven bijectively in two different ways in \cite[Theorem~6.1]{sortable} and \cite[Theorem~8.1]{sortable}.
In particular, the number of \mbox{$c$-sortable} elements is independent of~$c$.
Lemma~\ref{scs} gives a bijection between \mbox{$c$-sortable} elements shortened on the left by~$s$ and $scs$-sortable elements lengthened on the left by~$s$.
Thus there is a bijection between \mbox{$c$-sortable} elements~$w$ with $l(sw)>l(w)$ and $scs$-sortable elements~$x$ with $l(sx)<l(x)$.
Using Theorem~\ref{sublattice} and Lemma~\ref{s join 2}, one can show that $w\mapsto s\join w$ is such a bijection with inverse $x\mapsto (x)_{\br{s}}$.
\end{remark}

We conclude the section with an easy technical lemma which is useful in the proof of Theorem~\ref{camb}.

\begin{lemma}
\label{pidown same}
Let~$s$ be initial in~$c$ and let $y\in W$ have $l(sy)<l(y)$.
If $x\covered y$ but $x\neq sy$ then $\pidown^c(x)=\pidown^c(y)$ if and only if $\pidown^{scs}(sx)=\pidown^{scs}(sy)$.
\end{lemma}
\begin{proof}
Note that~$s$ is an inversion of~$y$ and the inversion sets of~$x$ and~$y$ differ by one reflection.
This reflection is not~$s$ because if so, we would have $x=sy$.
Thus $l(sx)<l(x)$.
So $\pidown^c(x)=s\cdot\pidown^{scs}(sx)$ and $\pidown^c(y)=s\cdot\pidown^{scs}(sy)$.
\end{proof}

\section{Lattice congruences of the weak order}
\label{cong sec}
In preparation for the proof (in Section~\ref{proof}) of Theorem~\ref{camb}, we review some well-known facts about congruences of finite lattices as well as some facts from~\cite{congruence} about the geometry underlying lattice congruences of the weak order on a finite Coxeter group.
(For more information about lattice congruences in a general setting, see~\cite{Gratzer}.)
We also discuss an observation (Observation~\ref{shard reflect}) which connects the geometry of lattice congruences to the action of a simple generator.
We remind the reader that the weak order on a Coxeter group~$W$ is a lattice if and only if~$W$ is finite.  
Thus the results and observations of this section should be applied only in the context of finite~$W\!$.\@

For a finite lattice $L$, let $\Con(L)$ be the set of congruences on $L$, partially ordered by refinement.
The poset $\Con(L)$ is a distributive lattice, 
and thus is completely specified by the subposet $\Irr(\Con(L))$ induced by the join-irreducible congruences.
(Recall that a join-irreducible element of a finite lattice is an element which covers exactly one element.) 
For a cover relation $x\covered y$ in $L$, a congruence is said to {\em contract} the edge $x\covered y$ if $x\equiv y$. For any cover relation $x\equiv y$ there is a unique smallest (i.e.\ finest in refinement order) congruence $\Cg(x\covered y)$ of $L$ which contracts that edge, and this congruence is join-irreducible 
in $\Con(L)$.
Every join-irreducible in $\Con(L)$ arises in this manner, and in fact, every join-irreducible congruence arises from a cover of the form $j_*\covered j$, where~$j$ is a join-irreducible element of $L$ and~$j_*$ is the unique element of $L$ covered by~$j$.
We say the congruence contracts~$j$ if it contracts the edge $j_*\covered j$ and write $\Cg(j)$ for the join-irreducible congruence $\Cg(j_*\covered j)$.
A congruence is determined by the set of join-irreducible elements it contracts.

The weak order on a finite Coxeter group~$W$ has a property called congruence uniformity or boundedness~\cite{bounded}, meaning in particular that the map $\Cg$ is a bijection from join-irreducibles of~$W$ to join-irreducibles of $\Con(W)$.
In what follows, we tacitly use the map $\Cg$ to blur the distinction between join-irreducible congruences and join-irreducible elements of the weak order on $W,$ so that $\Irr(\Con(W))$ is considered to be a partial order on the join-irreducibles of $W\!$.\@
To distinguish this partial order from the weak order on $W,$ we denote it ``$\le_{\Con}$.''
This partial order has the following description: 
$j_2\le_{\Con} j_1$ if and only if every congruence contracting $j_1$ must also contract $j_2$.
Thus congruences of~$W$ are identified with order ideals of contracted join-irreducibles in $\Irr(\Con(W))$.
(Order ideals in a poset are subsets~$I$ such that $x\in I$ and $y\le x$ implies $y\in I$.)

We now review a geometric characterization of the partial order $\le_{\Con}$ from~\cite{hyperplane}.
Continuing under the assumption that~$W$ is finite, we fix a reflection representation of~$W$ as a group of orthogonal transformations of a real Euclidean vector space~$V$ of dimension~$d$.
Each reflection $t\in T$ acts on~$V$ as an orthogonal reflection and the {\em Coxeter arrangement} for~$W$ is the collection~$\A$ of reflecting hyperplanes for reflections $t\in T$.
The hyperplanes in~$\A$ are permuted by the action of~$W\!$.\@
The {\em regions} of~$\A$ are the closures of the connected components of the complement \mbox{$V-(\cup\A)$} of~$\A$.
Choosing any region~$B$ to represent the identity element of $W,$ the elements of~$W$ are in one-to-one correspondence $w\mapsto wB$ with the regions of~$\A$.
The inversion set of an element is the {\em separating set} of the corresponding region~$R$
(the set of reflections whose hyperplanes separate~$R$ from~$B$). 
The weak order on~$W$ is thus containment order on separating sets.

Say a subset~$U$ of~$V$ is {\em above} a hyperplane~$H$ if every point in~$U$ is either contained in~$H$ or is separated from~$B$ by~$H\!$.\@
If~$U$ is above~$H$ and does not intersect~$H,$ then~$U$ is {\em strictly above}~$H\!$.\@
Similarly,~$U$ is {\em below}~$H$ if points in~$U$ are either contained in~$H$ or on the same side of~$H$ as~$B$, and~$U$ is strictly below~$H$ if it is disjoint from and below~$H\!$.\@

We say~$\A'$ is a {\em rank-two subarrangement} of~$\A$ if~$\A'$ consists of all the hyperplanes of~$\A$ containing some subspace of dimension $d-2$ 
and $|\A'|\ge 2$.
(The rank-two subarrangements of~$\A$ are exactly the collections of reflecting hyperplanes of the rank-two parabolic subgroups considered in~\cite{sortable}.)
For each rank-two subarrangement~$\A'$, there is a unique region $B'$ of~$\A'$ containing~$B$.
The two facet hyperplanes of $B'$ are called {\em basic} hyperplanes in~$\A'$.
(The basic hyperplanes of~$\A'$ are the reflecting hyperplanes for the canonical generators of the corresponding rank-two parabolic subgroups, as defined in~\cite[Section 1]{sortable}.)

We cut the hyperplanes of~$\A$ into pieces called {\em shards} as follows.
For each non-basic~$H$ in a rank-two subarrangement~$\A'$, cut~$H$ into connected components by removing the subspace $\cap\A'$ from~$H\!$.\@
Equivalently,~$H$ is cut along its intersection with either of the basic hyperplanes of~$\A'$. 
Do this cutting for each rank-two subarrangement, and call the closures of the resulting connected components of the hyperplanes {\em shards}.
(In some earlier papers \cite{hyperplane,hplanedim} where shards were considered, closures were not taken.) 
For illustrations of the shards for $W=A_3$ and $W=B_3$, see \cite[Figures~1 and~3]{congruence}.

The shards of~$\A$ are important because of their connection to the join-irreducible elements of $W\!$.\@ 
For any shard~$\Sigma$, let $U(\Sigma)$ be the set of {\em upper elements} for~$\Sigma$.
That is, $U(\Sigma)$ is the set of regions of~$\A$ having a facet contained in~$\Sigma$ such that the region adjacent through that facet is lower (necessarily by a cover) in the weak order.
We partially order $U(\Sigma)$ as an induced subposet of the weak order.
By \mbox{\cite[Proposition~2.2]{hplanedim}}, an element of~$W$ is join-irreducible in the weak order if and only if it is minimal in $U(\Sigma)$ for some shard~$\Sigma$.
Thus, to every join-irreducible~$j$ in~$W$ we associate a shard~$\Sigma_j$.
In fact, the map $j\mapsto\Sigma_j$ is a bijection \cite[Proposition 3.5]{congruence} between join-irreducibles of the weak order and shards in~$\A$.
The inverse map is written $\Sigma\mapsto j_\Sigma$.
We now use this bijection to describe $\Irr(\Con(W))$ in terms of shards.

For each shard~$\Sigma$, write~$H_\Sigma$ for the hyperplane in~$\A$ containing~$\Sigma$.
Define the {\em shard digraph} $\Sh(W)$ to be the directed graph whose vertices are the shards, and whose arrows are as follows:
Given shards~$\Sigma_1$ and~$\Sigma_2$, let~$\A'$ be the rank-two subarrangement containing~$H_{\Sigma_1}$ and~$H_{\Sigma_2}$.
There is a directed arrow $\Sigma_1\rightarrow \Sigma_2$ if and only if
\begin{enumerate}
\item[(i) ]$H_{\Sigma_1}$ is basic in~$\A'$ but~$H_{\Sigma_2}$ is not, and
\item[(ii) ]$\Sigma_1\cap\Sigma_2$ has dimension $d-2$.
\end{enumerate}
As explained in \cite[Section 8]{hyperplane}, this directed graph is acyclic.
The following is a special case of one of the assertions of \cite[Theorem 25]{hyperplane}.

\begin{theorem}
\label{shard}
The poset $\Irr(\Con(W))$ is isomorphic to the transitive closure of $\Sh(W)$ via the bijection $j\mapsto\Sigma_j$.
\end{theorem}

We interpret the transitive closure of $\Sh(W)$ as a poset by the rule that an arrow ``$\rightarrow$'' corresponds to an order relation ``$\ge$.'' 
In other words two join-irreducibles $j_1$ and $j_2$ in the weak order on~$W$ have $j_2\le_{\Con} j_1$ if and only if there is a directed path in $\Sh(W)$ from~$\Sigma_{j_1}$ to~$\Sigma_{j_2}$.

\begin{figure}[ht]

\scalebox{0.8}{\epsfbox{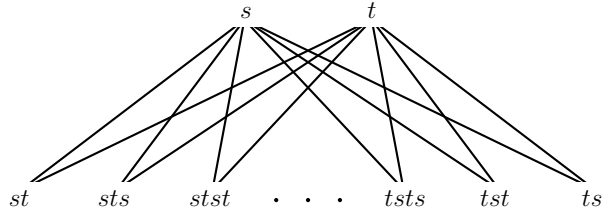}}

\caption{$\Irr(\Con(W))$ for~$W$ of rank 2}
\label{rk2}
\end{figure}

\begin{example}\rm
\label{Con rk 2}
For~$W$ of rank two and $S=\set{s,t}$ the partial order $\le_{\Con}$ on join-irreducibles is illustrated in Figure~\ref{rk2}.
(Cf. Figure~\ref{B2camb}.a in Section~\ref{intro}, which shows the weak order on the rank-two Coxeter group $B_2$.)
\end{example}

The following is a special case of \cite[Lemma~3.9]{congruence}.

\begin{lemma}
\label{source}
Let~$\Sigma$ be a shard.
The following are equivalent:
\begin{enumerate}
\item[(i) ]$\Sigma$ is a source in $\Sh(W)$.
\item[(ii) ]$\Sigma$ consists of an entire facet hyperplane of~$B$ (the region representing the identity element of~$W$).
\item[(iii) ]There is no facet of~$\Sigma$ intersecting the region for $j_\Sigma$ in dimension $d-2$.
\end{enumerate}
\end{lemma}

For $J\subseteq S$, let~$\A_J$ be the hyperplane arrangement associated to the standard parabolic subgroup $W_J$.
We think of the arrangement~$\A_J$ as a subset of the arrangement~$\A$ by letting $W_J$ inherit the reflection representation fixed for $W\!$.\@
Recall the notation $\br{s}=S-\set{s}$.
For a join-irreducible~$j$ in $W,$ write~$H_j$ for the hyperplane containing~$\Sigma_j$.
This is also the hyperplane separating~$j$ from~$j_*$, the unique element covered by~$j$.
The next two lemmas are special cases of \cite[Lemma 6.6]{congruence} and \cite[Lemma~6.8]{congruence} respectively.

\begin{lemma}
\label{parabolic lem}
For $s\in S$, suppose $H_1\in(\A-\A_{\br{s}})$ and $H_2\in\A_{\br{s}}$.
Let~$\A'$ be the rank-two subarrangement containing~$H_1$ and~$H_2$.
Then $(\A'\cap\A_{\br{s}})=\set{H_2}$ and~$H_2$ is basic in~$\A'$.
\end{lemma}

\begin{lemma}
\label{ji parabolic}
Let $s\in S$ and let~$j$ be a join-irreducible.
Then $H_j\in\A_{\br{s}}$ if and only if $j\in W_{\br{s}}$.
\end{lemma}

In what follows, note that each $s\in S$ is a join-irreducible element of~$W$ and~$H_s$ (the hyperplane separating~$s$ from the unique element $1$ covered by~$s$) happens to be the reflecting hyperplane for the reflection~$s$.
\begin{lemma}
\label{two piece}
Let $s\in S$ and suppose~$H$ is a hyperplane in~$\A$ such that $H\not\in\set{H_r:r\in S}$.
Then the following are equivalent:
\begin{enumerate}
\item[(i) ]$H$ is not basic in the rank-two subarrangement containing~$H$ and~$H_s$ but is basic in every other rank-two subarrangement containing~$H\!$.\@
\item[(ii) ] There are exactly two shards in~$H,$ and their intersection is $H\cap H_s$.
\item[(iii) ] $H\in\A_{\set{r,s}}$ for some $r\in S$ with $r\neq s$.
\end{enumerate}
\end{lemma}
\begin{proof}
The equivalence of (i) and (ii) is immediate from the definition of shards.

Suppose~$H$ satisfies (i) and (ii), but not (iii).
Let~$\Sigma$ be the shard in~$H$ below~$H_s$, let~$j$ be the unique minimal element of $U(\Sigma)$ and let $r\in S$ have $r\neq s$.
Since (iii) fails,~$H_r$,~$H_s$ and~$H$ are not in the same rank-two subarrangement.
Thus~$H_r$ intersects the interior of~$\Sigma$.
Since~$H$ is basic in the rank-two subarrangement it shares with~$H_r$, there is an element~$x$ of $U(\Sigma)$ such that the region for~$x$ is not above~$H_r$.
Therefore~$j$ is not above~$r$ in the weak order, for each $r\in S$ with $r\neq s$.
But also~$j$ is not above~$s$ because~$\Sigma$ is below the hyperplane~$H_s$.
The only element of the weak order not above some element of~$S$ is the identity $1$.
But $1$ is not join-irreducible, and this contradiction shows that (i) implies (iii).

Conversely, suppose~$H$ satisfies (iii).
Then the rank-two subarrangement containing~$H$ and~$H_s$ is $\A_{\set{r,s}}$.
By hypothesis,~$H$ is not basic in $\A_{\set{r,s}}$, because the basic hyperplanes of $\A_{\set{r,s}}$ are $H_r$ and $H_s$.
Every other rank-two parabolic~$\A'$ containing~$H$ also contains at least one hyperplane not in $\A_{\set{r,s}}$, so by Lemma~\ref{parabolic lem},~$H$ is basic in~$\A'$.
Thus $H$ satisfies (i).
\end{proof}

We conclude this section by discussing an observation that is crucial to the proof of Theorem~\ref{camb}.

\begin{observation}\rm
\label{shard reflect}
Let~$H$ be a hyperplane distinct from~$H_s$, let $sH\in\A$ be the hyperplane obtained from~$H$ by reflecting by~$s$, and let~$\A'$ be the rank-two subarrangement containing~$H,$ $sH$ and~$H_s$.
Then the decomposition of~$sH$ into shards agrees (via the reflection $s$) with the decomposition of~$H$ into shards except that one of~$H$ and~$sH$ may be cut by~$H_s$ while the other may not.
This exception occurs precisely when one of~$H$ and~$sH$ is basic in~$\A'$ but the other is not.  
(If both~$H$ and~$sH$ are basic in~$\A'$ then necessarily $H=sH$.)
\end{observation}

\begin{figure}[ht]
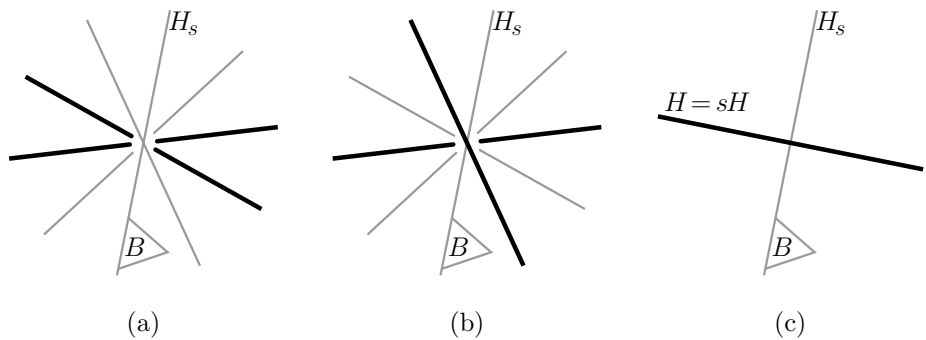

\centerline{
\begin{tabular}{ccccc}
\scalebox{.85}{\epsfbox{obs1a.ps}}&&\scalebox{.85}{\epsfbox{obs1b.ps}}&&\scalebox{.85}{\epsfbox{obs1c.ps}}\\[8 pt]
(a)&&(b)&&(c)
\end{tabular}
}
\caption{Illustrations for Observation~\ref{shard reflect}}
\label{obs1}
\end{figure}

Observation~\ref{shard reflect} is readily justified by the definition of shards.
Consider any rank-two subarrangement~$\A'$ containing~$H\!$.\@
The case where~$\A'$ also contains~$H_s$ is depicted in the three illustrations in Figure~\ref{obs1}.
Each illustration shows some of the hyperplanes in~$\A$ and indicates the position, with respect to these hyperplanes, of the region for~$B$.
We insert space between two shards in the same hyperplane to show where the cut is.
The hyperplanes~$H$ and~$sH$ are in black while the other hyperplanes are gray.
In the first two illustrations, the case $|\A'|=5$ is depicted.

If neither~$H$ nor~$sH$ is basic in~$\A'$ then both~$H$ and $sH$ are cut by~$H_s$ (possibly with $H=sH$) as exemplified in Figure~\ref{obs1}.a.
If exactly one of the two is basic in~$\A'$ then only the one that is non-basic is cut by~$H_s$, as exemplified in Figure~\ref{obs1}.b.
If both are basic then since~$H_s$ is also basic and~$\A'$ has only two basic hyperplanes, Figure~\ref{obs1}.c must apply.
These considerations explain how and why the shard decompositions of~$H$ and~$sH$ may or may not differ where these hyperplanes intersect~$H_s$.

\begin{figure}[ht]

\scalebox{0.76}{\epsfbox{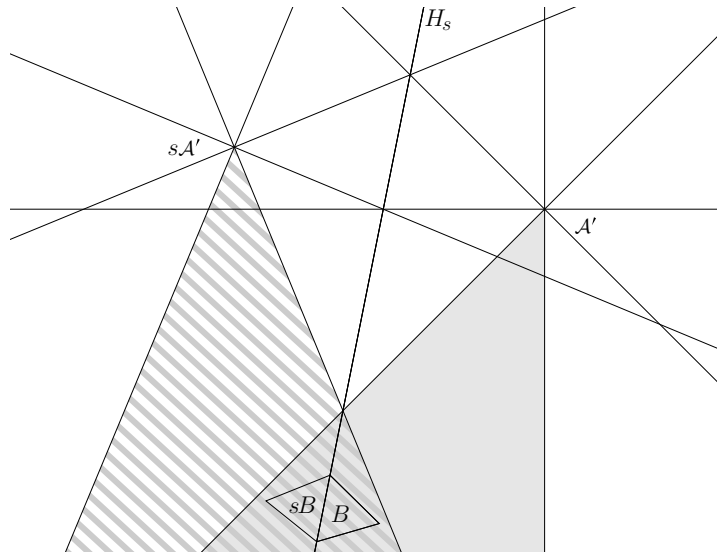}}

\caption{Another illustration for Observation~\ref{shard reflect}}
\label{obs2}
\end{figure}

To see why the shard decompositions of~$H$ and~$sH$ must agree away from~$H_s$, consider a rank-two subarrangement~$\A'$ containing~$H$  but not~$H_s$.
Let $s\A'$ be the rank-two subarrangement obtained from~$\A'$ by the reflection~$s$.
Since~$H_s$ is the only hyperplane separating the region~$B$  (representing $1$) from the region~$sB$ (representing~$s$) and $H_s\not\in\A'$,~$B$ and~$sB$ are contained in the same $\A'$-region $B'$.
The reflection~$s$ maps~$B$ to~$sB$, so~$B$ is contained in the $(s\A')$-region $sB'$.
Thus the reflection~$s$ maps the basic hyperplanes of~$\A'$ to the basic hyperplanes of $s\A'$.
(This is illustrated in Figure~\ref{obs2}, where $B'$ is shaded gray and $sB'$ is shaded in stripes.)
We conclude that~$H$ is non-basic in~$\A'$ (and is therefore cut by the basic hyperplanes of~$\A'$) if and only if~$sH$ is non-basic in $s\A'$ (and is therefore cut by the basic hyperplanes of $s\A'$).

\section{Cambrian lattices}
\label{camb sec}
In this section we define Cambrian congruences and Cambrian lattices and prove two lemmas which are useful in the proof of Theorem~\ref{camb}. 
We continue to restrict our attention to the case of a finite Coxeter group~$W\!$.\@

An {\em orientation} of the Coxeter diagram for~$W$ is obtained by replacing each edge of the diagram by a single directed edge, connecting the same pair of vertices in either direction.
Orientations of the Coxeter diagram correspond to Coxeter elements (cf.~\cite{shi}).
Specifically, for each pair of noncommuting simple generators~$s$ and $t$, the edge $s\,$---$\,t$ is oriented $s\!\to\! t$ if and only if~$s$ precedes $t$ in every reduced word for~$c$.
Each directed edge in the orientation corresponds to a pair of elements which are required to be congruent in the Cambrian congruence.
Specifically, if $s\!\to\! t$ then the requirement is that the element $t$ be congruent to the element with reduced word $tsts\cdots$ of length $m(s,t)-1$.
(Recall that $m(s,t)$ is the order of the product $st$ in~$W$.)
The {\em Cambrian congruence} associated to~$c$ is the smallest (i.e.\ finest as a partition) lattice congruence satisfying this requirement for each directed edge.
The Cambrian lattice associated to~$c$ is the quotient of the weak order on~$W$ modulo the Cambrian congruence.
For brevity, we refer to these as the $c$-Cambrian congruence and the $c$-Cambrian lattice.
Examples~\ref{B2lat} and~\ref{A3lat} illustrate these definitions.

The definition of the Cambrian congruence can be rephrased as follows:  
For each directed edge $s\!\to\! t$, we require that the join irreducibles $ts,tst,\ldots,tsts,...$ of lengths $2$ to $m(s,t)-1$ be contracted.
The join-irreducibles thus required to be contracted are called the {\em defining join-irreducibles} of the $c$-Cambrian congruence.
Requiring that these join-irreducibles be contracted specifies an order ideal of join-irreducibles under the partial order $\le_{\Con}$, namely the smallest order ideal containing the defining join-irreducibles.
This order ideal, interpreted as an order ideal in $\Irr(\Con(W))$, specifies the congruence.

\begin{example} \rm
\label{camb rk 2}
For~$W$ of rank two and $S=\set{s,t}$ the partial order $\le_{\Con}$ on join-irreducibles is illustrated in Figure~\ref{rk2} in Section~\ref{cong sec}.
The defining join-irreducibles for the Cambrian congruence for $c=st$ are the join-irreducibles of length greater than one whose unique reduced word starts with $t$.
This set of join-irreducibles is an order ideal.
Thus the Cambrian congruence has the interval $[t,tsts\cdots]$ as its only nontrivial congruence class, so that the bottom elements of the congruence classes are exactly the \mbox{$c$-sortable} elements (cf. Example~\ref{B2lat}).
\end{example}

We conclude this section with two lemmas which contribute to the proof of Theorem~\ref{camb}.

\begin{lemma}
\label{camb restrict}
Let $W_J$ be the standard parabolic subgroup generated by $J\subseteq S$ and let~$c'$ be the restriction of~$c$ to $W_J$.
Then $x\in W_J$ is contracted by the $c$-Cambrian congruence if and only if it is contracted by the $c'$-Cambrian congruence.
\end{lemma}
\begin{proof}
Recall that a lattice congruence is uniquely determined by the set of join-irreducibles it contracts.
As a special case of \cite[Lemma 6.12]{congruence}, the identity map embeds $\Irr(\Con(W_J))$ as an induced subposet of $\Irr(\Con(W))$, and the complement of this induced subposet is an order ideal.
The defining join-irreducibles of the $c'$-Cambrian congruence are exactly those defining join-irreducibles of the $c$-Cambrian congruence which happen to be in $W_J$.
Thus a join-irreducible in $W_J$ is below (in $\le_{\Con}$) some defining join-irreducible of the $c$-Cambrian congruence if and only if it is below some defining join-irreducible of the $c'$-Cambrian congruence.
\end{proof}

The {\em degree} of a join-irreducible~$j$ is the cardinality of the smallest $J\subseteq S$ such that $j\in W_J$.
By \cite[Lemma 6.12]{congruence}, if $j_2\le_{\Con} j_1$ in $\Irr(\Con(W))$ then the degree of $j_1$ is less than or equal to the degree of $j_2$.

\begin{lemma}
\label{camb para}
Let~$s$ be initial in~$c$.
If~$j$ is join-irreducible with $l(sj)>l(j)$ but $j\not\in W_{\br{s}}$ then~$j$ is contracted by the $c$-Cambrian congruence.
\end{lemma}
\begin{proof}
This is a modification of the proof of \cite[Theorem 6.9]{congruence}.
We argue by induction on the dual of $\Irr(\Con(W))$, the base case being where~$j$ is of degree one or two. 
The case where~$j$ is of degree one is vacuous.
If~$j$ is of degree two then~$j$ is a defining join-irreducible of the $c$-Cambrian lattice.

Let~$j$ be a join-irreducible of degree more than two satisfying the hypotheses of the lemma.
To accomplish the inductive proof, we need only find a join-irreducible~$j'$ above~$j$ in $\Irr(\Con(W))$ such that $l(sj')>l(j')$ but $j'\not\in W_{\br{s}}$.
Let~$H$ stand for~$H_j$ and let~$\Sigma$ be~$\Sigma_j$.
Since~$j$ is of degree more than two and each $s\in S$ has degree one,~$H$ is not a facet hyperplane of~$B$ (the region associated to $1$).

First, consider the case where~$H$ is basic in the rank-two subarrangement~$\A'$ containing~$H$ and~$H_s$.
By Lemma~\ref{source} it has a facet, as a polyhedral subset of~$H,$ and moreover, there is a facet of~$\Sigma$ intersecting the region for~$j$ in dimension $d-2$.
There are two hyperplanes in~$\A$ which define this facet, and since~$H_j$ is basic in~$\A'$, neither of these two hyperplanes is~$H_s$.
By Lemma~\ref{ji parabolic}, at least one of the two hyperplanes (call it $H'$) is not in~$\A_{\br{s}}$.
Some shard~$\Sigma'$ contained in~$H'$ arrows~$\Sigma$ in $\Sh(W)$ and thus intersects the region for~$j$ in dimension $d-2$.
Since the region for~$j$ is below~$H_s$, there is a region in $U(\Sigma')$ which is below~$H_s$.
In particular, $l(sj_{\Sigma'})>l(j_{\Sigma'})$, and since $H'\not\in\A_{\br{s}}$, Lemma~\ref{ji parabolic} says $j'\not \in W_{\br{s}}$.
Thus $j_{\Sigma'}$ is the desired~$j'$ in this case.

Next, consider the case where~$H$ is not basic in the rank-two subarrangement containing~$H$ and~$H_s$.
In particular,~$\Sigma$ is below~$H_s$.
By Lemmas~\ref{source} and~\ref{two piece} (since~$j$ has degree greater than two),~$\Sigma$ has a facet which is not defined by~$H_s$.
As in the previous case, this facet is defined by at least one hyperplane~$H'$ not in~$\A_{\br{s}}$.
Some shard~$\Sigma'$ in~$H'$ arrows~$\Sigma$, and since every region in $U(\Sigma)$ is below~$H_s$, the join-irreducible~$j'$ associated to~$\Sigma'$ has $l(sj')>l(j')$.
Finally,~$j'$ is not in $W_{\br{s}}$ by Lemma~\ref{ji parabolic}.
\end{proof}

\section{Sortable elements and Cambrian congruences}
\label{proof}
In this section, we prove Theorem~\ref{camb}, which states that the Cambrian congruence associated to~$c$ coincides with the congruence $\Theta_c$ of Theorem~\ref{cong} and that therefore the associated Cambrian lattice is the restriction of the weak order to \mbox{$c$-sortable} elements.
Part of the proof is accomplished by Theorem~\ref{cong}.
The strategy of the remainder of the proof is to relate the Cambrian congruences to the recursive structure of \mbox{$c$-sortable} elements as described in Lemmas~\ref{sc} and~\ref{scs}.

The definition of the Cambrian congruence can be further restated as follows:
The Cambrian congruence is the unique smallest congruence contracting all non-\mbox{$c$-sortable} join-irreducibles of degree 2.
(In fact, for any $s,t\in S$, every non-\mbox{$c$-sortable} element of $W_{\set{s,t}}$ is join-irreducible.)
The congruence $\Theta_c$ contracts all non-\mbox{$c$-sortable} join-irreducibles, and thus in particular contracts all non-\mbox{$c$-sortable} join-irreducibles of degree 2.
Therefore, by the definition of the Cambrian congruence, $\Theta_c$ is a weakly coarser congruence than the $c$-Cambrian congruence.
Since $\Theta_c$ does not contract any $c$-sortable join-irreducibles, we have the following corollary of Theorem~\ref{cong}:

\begin{cor}
\label{cong cor}
For a finite Coxeter group~$W$ and any Coxeter element~$c$, the $c$-Cambrian congruence does not contract any \mbox{$c$-sortable} join-irreducibles.
\end{cor}

In light of Corollary~\ref{cong cor} and the fact that a congruence is determined by the join-irreducibles it contracts, Theorem~\ref{camb} is equivalent to the following assertion.

\begin{prop}
\label{assertion}
For a finite Coxeter group~$W$ and any Coxeter element~$c$, the $c$-Cambrian congruence contracts every non-\mbox{$c$-sortable} join-irreducible.
\end{prop}

\begin{proof}
Let~$j$ be a non-\mbox{$c$-sortable} join-irreducible in~$W$ and fix an initial letter~$s$ of~$c$.
We argue by induction on the length of~$j$, on the rank of~$W$ and on the dual of the poset $\Irr(\Con(W))$ that~$j$ is contracted by the $c$-Cambrian congruence.

We first discuss the bases of the induction.
The proposition is trivial when $\rank(W)$ is $0$ or $1$, and is immediate from the definition for $\rank(W)=2$.
The base of the induction on length is the case $l(j)=2$.
The proposition holds in this case because any non-\mbox{$c$-sortable} element of length 2 is a defining join-irreducible for the $c$-Cambrian congruence.
The base of the induction on the dual of $\Irr(\Con(W))$ is the case where~$j$ is a defining join-irreducible (possibly of length greater than 2). 
In this case the proposition holds by definition.
Suppose now that~$j$ fits none of these base cases.

Consider the case where $l(sj)>l(j)$ and break into two subcases.
If $j\in W_{\br{s}}$ then by Lemma~\ref{sort para easy},~$j$ is not $sc$-sortable.
By induction on rank,~$j$ is contracted by the $sc$-Cambrian congruence, and by Lemma~\ref{camb restrict},~$j$ is contracted by the $c$-Cambrian congruence.
If $j\not\in W_{\br{s}}$ then by Lemma~\ref{camb para},~$j$ is contracted by the $c$-Cambrian congruence.

For the remaining case, $l(sj)<l(j)$, we apply Observation~\ref{shard reflect} and induction on length and on the dual of $\Irr(\Con(W))$.
To avoid repetitions of the phrase ``the region associated to,'' we identify each group element~$w$ with its corresponding region~$wB$.

Since~$s$ is \mbox{$c$-sortable} and~$j$ is not, $s\neq j$, so that $sj\neq 1$.
Therefore~$sj$ is join-irreducible in light of the isomorphism between $[s,w_0]$ and $[1,sw_0]$.
Lemma~\ref{scs} implies that $sj$ is not $scs$-sortable, and by induction on length,~$sj$ is contracted by the $scs$-Cambrian congruence.
Let~$\Sigma$ and~$H$ be the shard and hyperplane for~$sj$ and let~$\A'$ be the rank-two subarrangement containing~$H$ and~$H_s$.  
Since~$sj$ is below~$H_s$ and $sj\in U(\Sigma)$, the shard for~$\Sigma$ contains points strictly below~$H_s$.

First consider the case where~$\Sigma$ also contains a point strictly above~$H_s$.
We first claim that~$\A'$ contains at least one additional hyperplane.
To prove the claim, suppose to the contrary that $|\A'|=2$, as illustrated in Figure~\ref{impossible}.
In light of Observation~\ref{shard reflect} (particularly as illustrated in Figure~\ref{obs1}.c)~$\Sigma$ is fixed, as a set, by~$s$.
But then~$j$ is also in $U(\Sigma)$ so that by the uniqueness of join-irreducibles in $U(\Sigma)$, we must have $j=sj$, which is absurd.
This proves the claim.
\begin{figure}[ht]

\scalebox{0.8}{\epsfbox{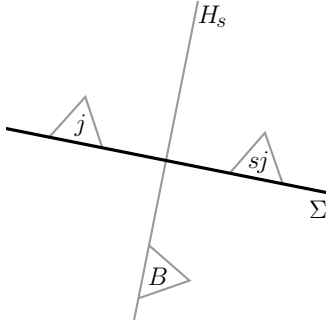}}

\caption{An impossible case in the proof of Proposition~\ref{assertion}}
\label{impossible}
\end{figure}

Let~$s\Sigma$ be the image of~$\Sigma$ under the reflection~$s$.
Since~$\Sigma$ contains both points above~$H_s$ and points below~$H_s$, the hyperplane~$H$ is basic in~$\A'$.
The other basic hyperplane is $H_s$, so in light of Observation~\ref{shard reflect} and the claim of the previous paragraph,~$s\Sigma$ is the union of two shards, one above~$H_s$ and one below.
Let~$\Sigma_1$ be the shard in~$s\Sigma$ below~$H_s$ and let~$\Sigma_2$ be the shard in~$s\Sigma$ above~$H_s$, as illustrated in Figure~\ref{above below}.

\begin{figure}[ht]

\scalebox{0.8}{\epsfbox{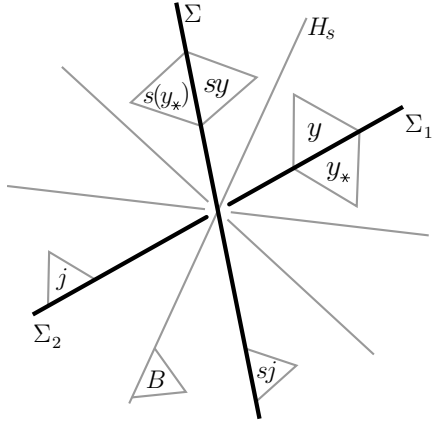}}

\caption{A case in the proof of Proposition~\ref{assertion}}
\label{above below}
\end{figure}

Now~$\Sigma$ arrows~$\Sigma_1$, so the join-irreducible~$y$ associated to~$\Sigma_1$ is also contracted by the $scs$-Cambrian congruence.
Thus by Corollary~\ref{cong cor},~$y$ is non-$scs$-sortable.
Since~$\Sigma_1$ is below~$H_s$, $l(sy)>l(y)$.
Let $y_*$ be the unique element covered by~$y$, so that~$\Sigma_1$ is the shard associated to the edge $y_*\covered y$.
Since~$y$ is a non-$scs$-sortable join-irreducible, it is contracted by $\Theta_{scs}$, or in other words $\pidown^{scs}(y_*)=\pidown^{scs}(y)$.
By Lemma~\ref{pidown same}, we have $\pidown^c(sy)=\pidown^c\left(s(y_*)\right)$, or in other words, the edge $s(y_*)\covered sy$ is contracted by $\Theta_c$.
But the shard associated to the edge $s(y_*)\covered sy$ is~$\Sigma$, so~$sj$ is contracted by $\Theta_c$.
By Corollary~\ref{cong cor},~$sj$ not \mbox{$c$-sortable}.
The shard associated to~$j$ is~$\Sigma_2$, which is arrowed to by~$\Sigma$.
By induction on length,~$sj$ is contracted by the $c$-Cambrian congruence, so~$j$ is also contracted by the $c$-Cambrian congruence.

By virtue of the cases we have already argued, we may now proceed under the assumptions that $l(sj)<l(j)$ and that the shard~$\Sigma$ associated to~$sj$ is below the hyperplane~$H_s$.

We claim that under these assumptions,~$sj$ is not a defining join-irreducible of the $scs$-Cambrian congruence.
Suppose to the contrary that~$sj$ is defining.
Then $sj\in W_{\br{s}}$, because all defining join-irreducibles for the $scs$-Cambrian congruence which are not in $W_{\br{s}}$ are shortened on the left by~$s$.
By Lemma~\ref{ji parabolic}, the hyperplane~$H$ containing~$\Sigma$ is in $\A_{\set{r_1,r_2}}$ with $s\not\in\set{r_1,r_2}$.
By Lemma~\ref{two piece},~$H$ is basic in the rank-two subarrangement containing~$H$ and~$H_s$ and~$\Sigma$ contains both points strictly above~$H_s$ and points strictly below~$H_s$.
This contradiction to the assumption proves the claim.

Thus our list of assumptions becomes:  $l(sj)<l(j)$;~$\Sigma$ is below~$H_s$; and~$sj$ is contracted by the $scs$-Cambrian congruence but is not a defining join-irreducible for the $scs$-Cambrian congruence.
By the last of these assumptions, there is a join-irreducible~$x$ which is contracted by the $scs$-Cambrian congruence and which arrows~$sj$ in $\Sh(W)$.
By Corollary~\ref{cong cor},~$x$ is non-$scs$-sortable.
Since~$x$ arrows~$sj$ in $\Sh(W)$, the corresponding shards~$\Sigma_x$ and~$\Sigma$ intersect in codimension~2, and the corresponding hyperplanes~$H_x$ and~$H$ are contained in a rank-two subarrangement~$\A'$ in which~$H_{x}$ is basic and~$H$ is not.
By the assumption that~$\Sigma$ is weakly below~$H_s$, there are two cases:  either $\Sigma_x\cap\Sigma$ is contained in~$H_s$ or $\Sigma_x\cap\Sigma$ contains a point strictly below~$H_s$.
We now complete the proof by considering these two cases.

\begin{figure}[ht]

\scalebox{0.8}{\epsfbox{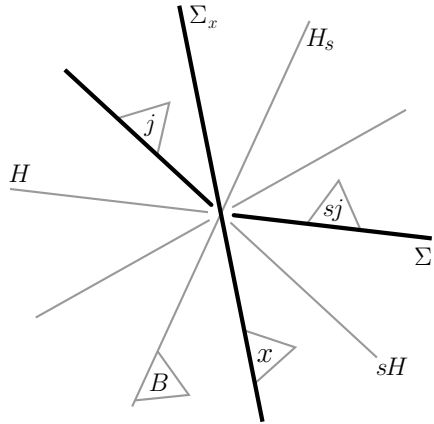}}

\caption{Another case in the proof of Proposition~\ref{assertion}}
\label{contained}
\end{figure}

The case where the intersection $\Sigma_x\cap\Sigma$ is contained in~$H_s$ is illustrated in Figure~\ref{contained}.
In this case,~$\A'$ contains~$H_s$, and since~$H_x$ is basic in~$\A'$, some upper region of~$\Sigma_x$ is below~$H_s$.
Thus $l(sx)>l(x)$.
Since~$x$ arrows~$sj$,~$\Sigma$ has a facet defined by~$H_x$, or equivalently by~$H_s$.
We claim that the hyperplane~$H$ is not the image of~$H_x$ under~$s$.
If it is, then by Observation~\ref{shard reflect}, the image~$s\Sigma$ of~$\Sigma$ under~$s$ is contained in a shard with points on both sides of~$H_s$ (namely, the shard~$\Sigma_x$).
However,~$s\Sigma$ is~$\Sigma_j$, which lies above~$H_s$ because $l(sj)<l(j)$.
This contradiction proves the claim.
Thus the hyperplane containing~$j$ is not basic in~$\A'$, so that~$j$ is arrowed to by~$x$, which is not \mbox{$c$-sortable}.

We next claim that~$x$ is not \mbox{$c$-sortable}.
If it is, then Lemma~\ref{sc} implies that~$x$ is an $sc$-sortable element of $W_{\br{s}}$.
But then by Lemma~\ref{sort para easy},~$x$ is $scs$-sortable as an element of~$W,$ contradicting what was established previously.
This contradiction establishes the claim.
By induction on the dual of $\Irr(\Con(W))$,~$x$ is contracted by the $c$-Cambrian congruence, and therefore so is~$j$.

Finally, we deal with the case where $\Sigma_x\cap\Sigma$ contains a point strictly below~$H_s$, as illustrated in Figure~\ref{not contained}.
In this case some upper region of~$\Sigma_x$ is below~$H_s$, so that the minimal upper region~$x$ has $l(sx)>l(x)$.
The reflection~$s$ relates~$\Sigma_x$ to a shard~$\Sigma'$ having $sx$ in its upper set.
The reflection~$s$ carries~$\Sigma$ into the shard~$\Sigma_j$.
Possibly~$\Sigma_x$ and~$\Sigma'$ differ to the extent allowed by Observation~\ref{shard reflect}, and similarly~$\Sigma$ and~$\Sigma_j$ may differ.
However, those differences occur only at the intersection with~$H_s$, in particular not affecting the following assertions.
\begin{figure}[ht]

\scalebox{0.8}{\epsfbox{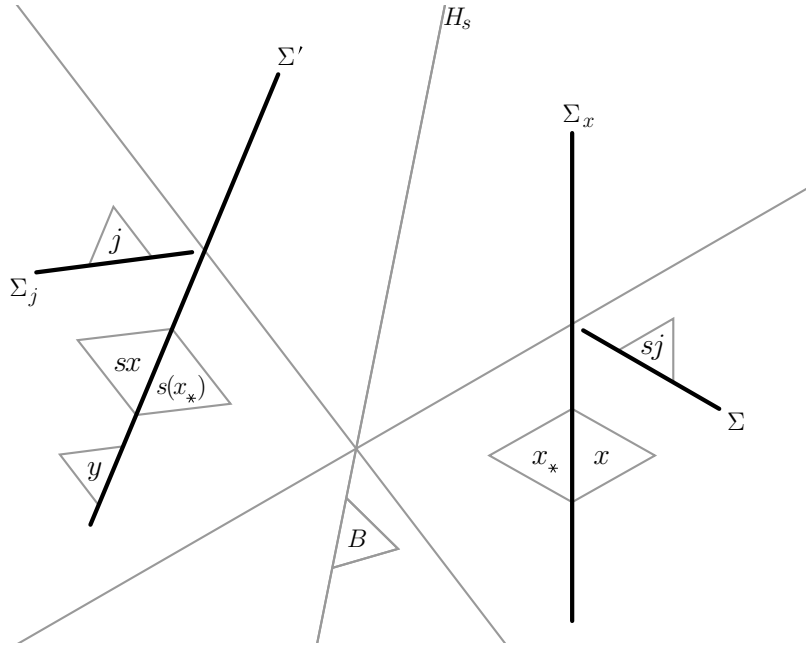}}

\caption{One last case in the proof of Proposition~\ref{assertion}}
\label{not contained}
\end{figure}

The intersection $\Sigma'\cap\Sigma_j$ has dimension $d-2$ because $\Sigma_x\cap\Sigma$ has dimension $d-2$.
As explained in connection with Observation~\ref{shard reflect}, reflection by~$s$ maps the basic hyperplanes of~$\A'$ to the basic hyperplanes of $s\A'$.
Thus~$\Sigma'$ arrows~$\Sigma_j$.
Let $x_*$ be the unique element covered by~$x$.
By Lemma~\ref{pidown same}, and the fact that $\Theta_{scs}$ contracts~$x$, the congruence $\Theta_c$ contracts the edge $sx_*\covered sx$.
Thus $\Theta_c$ contracts the join-irreducible~$y$ associated to~$\Sigma'$.
By Corollary~\ref{cong cor}, the join-irreducible~$y$ is non-\mbox{$c$-sortable} and by construction~$y$ is above~$j$ in $\Irr(\Con(W))$.
By induction on the dual of $\Irr(\Con(W))$,~$y$ is contracted by the $c$-Cambrian congruence, and therefore so is~$j$.
\end{proof}

This concludes the proof of Theorem~\ref{camb}.

\section*{Acknowledgments}
I am grateful to David Speyer, John Stembridge and Hugh Thomas for helpful conversations.
I particularly thank David Speyer for pointing out, in the course of joint work on~\cite{camb_fan} that Lemma~\ref{cov w br s} was not immediately obvious, and for cooperating in writing down the proof.

\end{document}